\documentclass[12pt,a4paper]{article}

\usepackage{amssymb}

\newtheorem{theorem}{\bf Theorem}[section]
\newtheorem{algorithm}{\bf Algorithm}[section]
\newtheorem{proposition}{\bf Proposition}[section]

\title{Random Permutations, Random Sudoku Matrices  and Randomized Algorithms}
\author{Krasimir Yordzhev}
\date{}

\begin{document}
\maketitle
\begin{center} {\em
Faculty of Mathematics and Natural Sciences\\
South-West University, Blagoevgrad, Bulgaria} \\
E-mail: yordzhev@swu.bg
\end{center}

\begin{abstract}Some randomized algorithms, used to obtain a random   $n^2 \times n^2$ Sudoku matrix, where $n$ is a natural number, is reviewed in this study. Below is described the set $\Pi_n$ of all $(2n) \times n$ matrices, consisting of elements of the set $\mathbb{Z}_n =\{ 1,2,\ldots ,n\}$, such that every row is a permutation. It is proved that such  matrices would be particularly useful in developing efficient algorithms in generating Sudoku matrices.
An algorithm to obtain random $\Pi_n$ matrices  is presented in this paper.  The algorithms are evaluated according to two criteria - probability evaluation, and time evaluation. This type of criteria is interesting from both theoretical and practical
point of view because they are particularly useful in the analysis of computer programs.

\end{abstract}



\section{Introduction}\label{intr}
This work can be particularly useful to future computer engineers and their lecturers (see also \cite{Ishrat,SRINIVASAN}).

Let $\mathfrak{M}$ be a finite set. \emph{A Random objects generator of $\mathfrak{M}$}  is every algorithm $\mathcal{A}_\mathfrak M$  randomly generating any element of $\mathfrak M$, while elements generated by a random objects generator will be called \emph{random elements } of $\mathfrak M$, i.e. random numbers, random matrices, random permutations, etc. We take for granted that probabilities to obtain different random elements of  $\mathfrak M$ by means of $\mathcal{A}_\mathfrak M$ are equal, and are also equal to $\displaystyle \frac{1}{\left| \mathfrak{M} \right|}$. We denote the time that the random objects generator needs to obtain a random element of $\mathfrak M$ with $T(\mathcal{A}_\mathfrak{M} )$.

A \emph{randomized algorithm} is an algorithm which employs a degree of randomness as part of its logic. The randomized algorithms are very useful mathematical methods for solving a class of problems, which uses a random objects generator \cite{gentle1,gentle2}.
In computing, a \emph{Monte Carlo algorithm} is a randomized algorithm whose running time is deterministic, but whose output may be incorrect with a certain probability.
The related class of \emph{Las Vegas algorithms} is also randomized, but in a different way: they take an amount of time that varies randomly, but always produce the correct answer.
A Monte Carlo algorithm runs for a fixed number of steps, and produces an answer that is correct with probability. A Las Vegas algorithm always produces the correct answer and its running time is a random variable.

For the purpose of this study, some randomized algorithms are used to solve the following class of problems: Let $n$ and $m=m(n)$ be natural numbers. Let us take the set $\mathcal{U}$, consisting of objects, dependent on $m$ parameters, where every parameter belongs to the finite set $\mathfrak M$. Let $\mathcal{V} \subset \mathcal{U}$. The problem is to obtain (at least one) object, which belongs to the set $\mathcal{V}$. The number of the elements of the sets $\mathcal{U}$ and $\mathcal{V}$ depends only on the parameter $m$, which is an integer function of the argument $n$. The standard that is most often used to illustrate a randomized algorithm (in particular a Las Vegas algorithm) develops as follows:

\begin{algorithm}\label{alg1}   \hfill

1) We obtain consequently $m=m(n)$ random elements of $\mathfrak M$ using random objects generator $\mathcal{A}_\mathfrak M$ with the help of which we initialize parameters of the object $u\in\mathcal{U}$;

2) We check if $u\in\mathcal{V}$. If the answer is no, everything is repeated.
\end{algorithm}

In other words, if we already have a random objects generator, a randomized algorithm can be used as a generator of more complex random objects. The benefit of these algorithms is that they provide a clear and simple description of any particular algorithm. Thus the randomized algorithms could be the basis for the development of different algorithms, which solve one and the same problem, and whose efficiency may vary in any particular case, as it is shown below.

The efficiency of Algorithm \ref{alg1} depends on the particular case in which it is used and can be evaluated according to the following criteria \cite{petrov}:

\textbf{Probability evaluation}: If $p(n)$ denotes the probability after generating $m=m(n)$ random elements of $\mathfrak{M}$ of obtaining an object of $\mathcal{V}$, then according to the classical probability formula:
\begin{equation}\label{p(n)}
p(n) = \frac{\left| \mathcal{V}\right|}{\left| \mathcal{U}\right|}
\end{equation}

\textbf{Time evaluation}: We denote by $\tau (n)$ the time needed to execute one iteration of Algorithm \ref{alg1}. Then
\begin{equation}\label{tau}
\tau (n) =m(n) T(\mathcal{A}_\mathfrak{M})+\theta (n) ,
\end{equation}
where $\theta (n)$ is the time to examine if the obtained object belongs to the set $ \mathcal{V}$.

It is obvious that the efficiency of Algorithm \ref{alg1} will be directly proportional to $p(n)$ and inversely proportional to $\tau (n)$.

The cases in which probability evaluation is equal to 1, i.e. the cases in which the algorithm is constructed directly to obtain element of the set $\mathcal V$ and there is no need of belonging examination, are of great interest, as only one iteration is implemented then, i.e. there is no repetition. The algorithm for  obtaining  random permutations, based on a randomized algorithm with a probability evaluation is equal to  1, which is more efficient than the other algorithm with probability evaluation $\displaystyle p (n) = \frac{n!}{n^n}$, is described in section \ref{sec2}.

Let $k$ be an integer. We denote by $\mathbb{Z}_k$ the set of the integers
$\mathbb{Z}_k =\left\{ 1,2,\ldots ,k \right\}$ and by $\mathcal{S}_n$ the set of all permutations of elements of $\mathbb{Z}_n$.  As it is well known
$|\mathcal{S}_n |=n! $ .

There are standard procedures for obtaining random numbers of the set $\mathbb{Z}_k$ in most of the programming environments. We take this statement for granted and we will use it in the following examinations. Let $\mathcal{A}_k$ be a similar procedure. We consider that $T (\mathcal{A}_k )\approx T (\mathcal{A}_l) \approx t_0 = \textrm{Const}$ for $k\ne l$ in the current study.

Let $P_{ij}$, $1\le i,j\le n$ are $n^2$ in number square $n\times n$ matrices, whose elements belong to the set $\mathbb{Z}_{n^2} =\{ 1,2,\ldots ,n^2 \}$. Then $n^2 \times n^2$ matrix $P =\left[ P_{ij} \right]$ is called a \emph{Sudoku Matrix}, if every row, every column and every submatrix $P_{ij}$, $1\le i,j\le n$ make permutation of the elements of set $\mathbb{Z}_{n^2}$, i.e. every number $s\in \{ 1,2,\ldots ,n^2 \}$ is present only once in every row, every column and every submatrix $P_{ij}$. Submatrices $P_{ij}$ are called blocks of $P$.

Below we will illustrate the above mentioned ideas by analyzing an arbitrary permutation of $n$ elements, an arbitrary $n^2 \times n^2$ Sudoku matrix  and an arbitrary $(2n) \times n$ matrix  with $2n$ rows and $n$ columns, every column of which is a permutation of $n$ elements, which are obtained by randomized algorithms.

We will prove that the problem for obtaining ordered $n^2$ - tuple of $(2n) \times n$ matrices, every row of which is a permutation of elements of $\mathbb{Z}_n$ is equivalent to the problem of generating a Sudoku matrix. We will analyze some possible algorithms for generating a random Sudoku matrix.

The randomized algorithms are very often used to solve problems, which are proved to be NP-complete. For detailed information about NP-complete problems and their application see \cite{garey} or \cite{hopcroft}. A proof that a popular Sudoku puzzle is NP-complete is given in \cite{yato} and \cite{yatoseta}. How to create computer program for Sudoku solving, using the concept set combined with the trial and error method is described in \cite{yorkost}.

\section{Random permutations}\label{sec2}
\par If there is a permutation of all elements of the set $\mathbb{Z}_n =\{ 1,2,\ldots ,n\}$ then $m(n)=n$.

We denote by $p_1 (n)$ the probability to obtain a random permutation of $\mathcal{S}_n$ with the help of Algorithm \ref{alg1}. Then according to formula  (\ref{p(n)}) we obtain:

\begin{equation}\label{p1}
p_1 (n) = \frac{n!}{n^n}
\end{equation}

\begin{proposition}\label{prop1}
There is an algorithm working in time $O(n)$ and checking if ordered $n$-tuple $\rho =\langle a_1 ,a_2 ,\ldots a_n \rangle$, $a_i \in \mathbb{Z}_n$ is a permutation, where $i=1,2,\ldots ,n$.
\end{proposition}

 To prove this we shall use the following algorithm, which obviously works in time $O(n)$:

\begin{algorithm}\label{alg2}
Check if given $n$-tuple $\rho =\langle a_1 ,a_2 ,\ldots a_n \rangle$, $a_i \in \mathbb{Z}_n$ is a permutation.\\
1) We declare array of $n$ elements
$v[1] ,v[2] ,\ldots ,v[n]$
and initialize all of its elements with 0;\\
2) \textbf{for} $i=1,2,\ldots ,n$ \textbf{do}

\textbf{begin}\\
3) $v[a_i ] := v[a_i ] +1$;\\
4) \textbf{if} $v[a_i ] >1$ \textbf{then}  $\{ \langle a_1 ,a_2 ,\ldots ,a_i ,\ldots ,a_n \rangle \notin \mathcal{S}_n ;\;   \textbf{exit}; \} \quad$  /* because the number $a_i$ is found more than once in $\rho$ and exit of the algorithm with the negative output */

\textbf{end}
\end{algorithm}
\hfill $\Box$

Let $\tau_1 (n)$ denote the time for implementing one iteration of Algorithm \ref{alg1} when a random permutation of elements of $\mathbb{Z}_n$ is obtained, and let $\theta_1 (n) $ denote the time for checking whether an arbitrary $n$-tuples of numbers of $\mathbb{Z}_n $ belongs to $\mathcal{S}_n$. Then, having in mind the formula (\ref{tau}) and Proposition \ref{prop1} we obtain the following time evaluation:
\begin{equation}\label{tau1}
\tau_1 (n) =n T(\mathcal{A}_n)+\theta_1 (n) =nt_0 +O(n)=O(n).
\end{equation}

The following algorithm is also randomized  (random numbers are generated), but its probability evaluation is equal to 1, i.e. in Algorithm \ref{alg1} step 2 is not implemented, because when the first random $n$ numbers are generated the obtained ordered $n$-tuple is a permutation.

\begin{algorithm}\label{alg3} Obtaining random permutation $\rho =\langle a_1 ,a_2 ,\ldots ,a_i ,\ldots ,a_n \rangle \in \mathcal{S}_n ,$ where $a_i \in \mathbb{Z}_n$, $i=1,2,\ldots ,n$, $a_i \ne a_j$ when $i\ne j$.\\
1) We declare array with $n$ elements
$v[1] ,v[2] ,\ldots ,v[n]$;\\
2) \textbf{for} $k=1,2,\ldots ,n $ \textbf{do} $v[k]:=k$;\\
3) \textbf{for} $k=1,2,\ldots ,n $ \textbf{do}

\textbf{begin}\\
4) We generate random number $x\in \mathbb{Z}_{n-k+1} =\{ 1,2,\ldots ,n-k+1\}$;\\
5)  $ a_k := v[x] ;$\\
6) \textbf{for} $j=x,x+1,\ldots ,n-k$ \textbf{do} $v[j]:=v[j+1] ; \quad$ /* We delete the element $v[x]$ and reduce the number of the elements of the array with 1*/

\textbf{end}
\end{algorithm}

It is obvious that the following proposition is true:

\begin{proposition}\label{prop2}
Algorithm \ref{alg3} which obtains random permutation has \textbf{probability evaluation}:
\begin{equation}\label{p2}
p_2 (n) =1
\end{equation}
and \textbf{time evaluation}:
\begin{equation}
\tau_2 (n) =t_0 \left[ O(n) +O(n-1)+\cdots +O(1)\right] =O(n^2 ) .
\end{equation}
\hfill $\Box$
\end{proposition}

 The probability evaluation of Algorithm \ref{alg3} is equal to 1.  The  probability evaluation of Algorithm \ref{alg1} is equal to $\displaystyle p_1 (n) = \frac{n!}{n^n} < 1$ for $n\ge 2$ and $\displaystyle \lim_{n\to\infty} p_1 (n) =0$. Therefore Algorithm \ref{alg3} is more efficient in obtaining random permutations than Algorithm \ref{alg1} with respect to the probability  evaluation.

\section{Random $(2n)\times n$ matrices, every row of which is a permutation of elements of $\mathbb{Z}_n$ }\label{sec3}

Let $\Pi_n$ denote the set of all $(2n)\times n$ matrices, which are also called $\Pi_n$ matrices, in which every row is a permutation of all elements of ${Z}_n$. In this case $\mathfrak{M}=\mathbb{Z}_n$ and $m(n)=2n^2$. It is obvious that
\begin{equation}\label{Pi_n}
\left| \Pi_n \right| =\left( n! \right)^{2n}
\end{equation}

It is easy to see that the following proposition is true:

\begin{proposition}\label{Pin_alg1}
When we obtain random $\Pi_n$ matrix with the help of Algorithm \ref{alg1} the following evaluations can be observed:

{\bf Probability evaluation:}
\begin{equation}\label{p3(n)}
p_3 (n) = \frac{\left| \mathcal{V}\right|}{\left| \mathcal{U}\right|} =\frac{(n! )^{2n}}{n^{2n^2}}
\end{equation}

{\bf Time evaluation:}
\begin{equation}\label{tau3}
\tau_3 (n) =m(n) T(\mathcal{A}_n)+2n\tau_1 (n)=2n^2 t_0 +2nO(n)=O(n^2 ) ,
\end{equation}
where $\tau_1 (n)$ is obtained according to formula (\ref{tau1}).

\hfill $\Box$
\end{proposition}

As it could be seen below (Proposition \ref{prop4}) the following algorithm is more efficient than Algorithm \ref{alg1} in obtaining random $\Pi_n$ matrix according to the probability evaluation.

\begin{algorithm}\label{alg4}
Obtaining random $\Pi_n$ matrix with probability evaluation equal to 1.

1) \textbf{for} $k=1,2,\ldots ,n , n+1,\ldots 2n$ \textbf{do} We obtain random permutation with the help of Algorithm \ref{alg3}, which will be the $k$-th row of the matrix;
\end{algorithm}

Practically, Algorithm \ref{alg4} repeats Algorithm \ref{alg3} $2n$ times. As it is stated in Proposition \ref{prop2} that Algorithm \ref{alg3} has a probability evaluation equal to 1, it logically follows, that Algorithm \ref{alg4} will have probability evaluation equal to 1. Then, we obtain the following proposition could be obtained:

\begin{proposition}\label{prop4}
Algorithm \ref{alg4} which obtains random $\Pi_n$ matrix has \textbf{probability evaluation}:
\begin{equation}\label{p4}
p_4 (n) =1
\end{equation}
and \textbf{time evaluation}:
\begin{equation}\label{tau4}
\tau_4 (n) =2nO(n^2 ) =O(n^3 ) .
\end{equation}
\hfill $\Box$

\end{proposition}

As we can see below $\Pi_n$ matrices can successfully be used to create algorithms that are efficient in developing Sudoku matrices.

A matrix, whose elements are equal to 0 or 1 is called \emph{binary}. A square binary matrix is called \emph{ permutation}, if there is only one 1 in every row and every column of the matrix. Let $\Sigma_{n^2}$ denote the set of all permutation $n^2 \times n^2$ matrices of the following type
\begin{equation}\label{matrA}
A =
\left[
\begin{array}{cccc}
A_{11} & A_{12} & \cdots & A_{1n} \\
A_{21} & A_{22} & \cdots & A_{2n} \\
\vdots & \vdots & \ddots & \vdots \\
A_{n1} & A_{n2} & \cdots & A_{nn}
\end{array}
\right] ,
\end{equation}
where for every $s,t\in \{ 1,2,\ldots ,n\}$  $A_{st}$ is a square $n\times n$ binary submatrix (block) with only one element equal to 1, and the rest of the elements are equal to 0.
As it is proved in \cite{dahl}
\begin{equation}\label{sssssss}
\left| \Sigma_{n^2} \right| =\left( n! \right)^{2n}
\end{equation}

Therefore, if a random $\Sigma_{n^2}$ matrix obtained by means of Algorithm \ref{alg1} the following probability evaluation is valid:

\begin{equation}\label{ppp5}
p_5 (n)=\frac{(n!)^{2n}}{2^{(n^2 )^2}} =\frac{(n!)^{2n}}{2^{n^4}}
\end{equation}

In order to obtain a random $\Sigma_{n^2}$ matrix, we have to generate $m(n)=(n^2 )^2 =n^4$ random numbers, which belong to the set $\{ 0,1\}$. Hence, whatever randomized algorithm  is used, the result is the following time evaluation:
\begin{equation}\label{ttt5}
\tau_5 (n) =m(n)T(\mathcal{A}_2) +\theta_5 (n)=n^4 t_0 +\theta_5 (n) =O(n^z ),\quad z\ge 4,
\end{equation}
where $\theta_5 (n)$ is check-up time, if the given binary matrix belongs to the set $\Sigma_{n^2}$.

In order to check if the given binary matrix is $\Sigma_{n^2}$ we can use, for example, the following algorithm, working in time $O(n^4 )$, i.e. $\theta_5  =O(n^4 )$ and therefore $z=4$.

\begin{algorithm}\label{Sigman^2}
Check if binary $n^2 \times n^2$ matrix $B=(b_{ij} )\in \Sigma_{n^2}$.\\
1) \textbf{for} $i=1,2,\ldots ,n^2$ \textbf{do}

\textbf{begin}\\
2) r:=0;\\
3) c:=0;\\
4) \textbf{for} $j=1,2,\ldots ,n^2$  \textbf{do}

 \hspace{1cm} \textbf{begin}\\
5) \hspace{1cm} $r:=r+b_{ij} $;\\
6) \hspace{1cm} \textbf{if}  $r>1$ \textbf{then} $B$ is not a permutation matrix and exit of the algorithm;\\
7) \hspace{1cm} $c:=c+b_{ji}$;\\
8 \hspace{1cm} \textbf{if} $c>1$ \textbf{then} $B$ is not a permutation matrix and exit of the algorithm;

 \hspace{1cm} \textbf{end};\\
9) \textbf{if} $r=0$ \textbf{or} $c=0$ \textbf{then} $B$ is not a permutation matrix and exit of the algorithm;

\textbf{end};\\
10) \textbf{for} $s=0,1,\ldots ,n-1$ \textbf{do}\\
11) \textbf{for} $t=0,1,\ldots ,n-1$ \textbf{do}

\textbf{begin}\\
12) $x:=0$;\\
13) \textbf{for} $i=1,2,\ldots ,n$ \textbf{do}\\
14) \textbf{for} $j=1,2,\ldots ,n$ \textbf{do} $x:=x+b_{sn+i\; tn+j}$;\\
15) \textbf{if} $x\ne 1$ \textbf{then} $B\notin \Sigma_{n^2}$  and exit of the algorithm;

\textbf{end}.
\end{algorithm}

When we compare (\ref{ppp5}) with (\ref{p3(n)}) and (\ref{p4}), as well as (\ref{ttt5}) with (\ref{tau3}) and (\ref{tau4}) we may assume that algorithms which use random $\Pi_n$ matrices are expected to be more efficient regard to probability and time evaluation than algorithms using random $\Sigma_{n^2}$ matrices to solve similar problems. This gives grounds for further examination of the $\Pi_n$ matrices' properties.

Two $\Sigma_{n^2}$  matrices $A=(a_{ij} )$ and $B=( b_{ij} )$, $1\le i,j\le n^2$ are called \emph{disjoint}, if there are no integers $i,j\in \{ 1,2,\ldots ,n^2 \}$  such that $a_{ij} =b_{ij} =1$.

We will give a little bit more complex definition of the term ''disjoint'' regarding $\Pi_n$ matrices. Let $C=(c_{ij} )$ and $D=(d_{ij} )$, $1\le i\le 2n$, $1\le j\le n$ be two $ \Pi_n$ matrices. We regard $C$ and $D$ as \emph{disjoint} matrices, if there are no natural numbers $s,t\in \{ 1,2,\ldots n\}$ such that the ordered pair $\langle c_{st} ,c_{n+t\; s} \rangle$ is equal to the ordered pair $\langle d_{st} ,d_{n+t\; s} \rangle$.

The following obvious proposition is given in \cite{dahl}:

\begin{proposition}\label{disj} \cite{dahl}
A square $n^2 \times n^2$ matrix $P$ with elements of  $\mathbb{Z}_{n^2} =\{ 1,2,\ldots ,n^2 \}$ is Sudoku matrix if and only if there are mutually disjoint matrices $A_1 ,A_2 ,\ldots ,A_{n^2} \in\Sigma_{n^2}$ such that $P$ can be presented as follows:
$$P=1\cdot A_1 +2\cdot A_2 +\cdots +n^2 \cdot A_{n^2}$$
\hfill $\Box$
\end{proposition}

The relationship between $\Pi_n$ matrices and Sudoku matrices is illustrated by the following theorem, considering Proposition \ref{disj}.

\begin{theorem}\label{th1}
There is a bijective map from $\Pi_n$ to $\Sigma_{n^2}$ and the pair of disjoint matrices of $\Pi_n$ corresponds to the pair of disjoint matrices of $\Sigma_{n^2}$
\end{theorem}

Proof. Let $P=(p_{ij} )_{2n\times n} \in \Pi_n$. We obtain an unique matrix of $\Sigma_{n^2}$ from $P$ by means of the following algorithm:

\begin{algorithm}\label{alg5}
Obtain just one  $\Sigma_{n^2}$ matrix, if $P=(p_{ij} )_{2n\times n} \in \Pi_n$ is given .\\
1) \textbf{for}  $ s = 1,2,\ldots ,n$ \textbf{do}\\
2) \textbf{for}  $ t = 1,2,\ldots ,n$ \textbf{do}

\textbf{begin}\\
3) $k:=p_{st}$;\\
4) $l:=p_{n+t\; s}$;\\
5) We obtain $n\times n$ matrix $A_{st} =(a_{ij})_{n\times n}$ such that $a_{kl}= 1$ and $a_{ij} =0$ in all other occasions;

\textbf{end};\\
6) We obtain matrix $A$ according to formula (\ref{matrA});
\end{algorithm}

Let $s\in \mathbb{Z}_n =\{ 1,2,\ldots ,n\}$. Since ordered $n$-tuple $\langle p_{s1} ,p_{s2} ,\ldots ,p_{sn} \rangle$ which is $s$-th row of the matrix $P$ is a permutation, then in every row of $n\times n^2$ matrix
$$
R_s =
\left[
\begin{array}{cccc}
A_{s1} & A_{s2} & \cdots & A_{sn}
\end{array}
\right]
$$
there is only one 1. For every $j=1,2,\ldots ,n$ $A_{sj}$ is binary $n\times n$ matrix in this case.

Similarly,  since ordered $n$-tuple $\langle p_{n+t\; 1} ,p_{n+t\; 2} ,\ldots ,p_{n+t\; n} \rangle$ which is $(n+t)$-th row of $P$ is a permutation for every $t\in \mathbb{Z}_n$, then in every column of $n^2 \times n$ matrix
$$
C_t =
\left[
\begin{array}{c}
A_{1t}  \\
A_{2t} \\
\vdots \\
A_{nt}
\end{array}
\right]
$$
there is only one 1, where $A_{it}$, $i=1,2,\ldots ,n$ is a binary $n\times n$ matrix. Hence, the matrix $A$ which is obtained with the help of Algorithm \ref{alg5} is $\Sigma_{n^2}$ matrix.

Since a unique matrix  of $\Sigma_{n^2}$ is obtained for every $P\in \Pi_n$ by means of Algorithm \ref{alg5}, then this algorithm provides a description of the map $\varphi \; :\; \Pi_n \to \Sigma_{n^2}$. It is easy to see that if there are given different elements of $\Pi_n$, we can use Algorithm \ref{alg5} to obtain different elements of $\Sigma_{n^2}$. Hence, $\varphi$ is an injection. But according to formulas (\ref{Pi_n}) and (\ref{sssssss}) $\left| \Sigma_{n^2} \right| =\left| \Pi_n \right|$, whereby it follows that $\varphi$ is a bijection.

Analyzing Algorithm \ref{alg5}, we arrived at the conclusion that $P$ and $Q$ are disjoint matrices of $\Pi_n$ if and only if $\varphi (P)$ and $\varphi (Q)$ are disjoint matrices of $\Sigma_{n^2}$ according to the above mentioned definitions. The theorem is proved.

\hfill $\Box$

\section{Conclusion}\label{sec4}

Let  $\mathcal{M}_{n^2} $ be the set of all square $n^2 \times n^2$ matrices with elements of the set $\mathbb{Z}_{n^2} =\{ 1,2,\ldots ,n^2 \}$, and let $\sigma_{n}$ be the number of all $n^2 \times n^2 $ Sudoku matrices. Obviously, $\left| \mathcal{M}_{n^2} \right| =(n^2 )^{n^2} =n^{2n^2}$. Then if we use Algorithm \ref{alg1} to obtain random Sudoku matrix. According to formula (\ref{p(n)}) there is the following probability evaluation:
\begin{equation}\label{p6}
p_6 (n) =\frac{\sigma_n}{n^{2n^2}}
\end{equation}

When $n=2$  $\sigma_2 = 288$ \cite{yorkost}.
When $n=3$, there are exactly
$\sigma_3 = 6\; 670\; 903\; 752\; 021\; 072\; 936\; 960 \approx 6.671\times 10^{21}$ in number Sudoku matrices \cite{Felgenhauer}.
As far as the author of this study knows, there is not a universal formula for the number $\sigma_n$ of  Sudoku matrices with every natural number $n$. We consider it as an open problem in mathematics.

If we employ random methods to create the matrix $P\in \mathcal{M}_{n^2}$ with elements of $\mathbb{Z}_{n^2}$, then according to Algorithm \ref{alg1} we need to verify if every row, every column and every block of $P$ is a permutation of elements of $\mathbb{Z}_{n^2}$. According to Proposition \ref{prop1}, every verification can be done in $O(n)$ time. Hence, when we employ Algorithm \ref{alg1} to obtain a random Sudoku matrix we will obtain the following time evaluation emerges:
\begin{equation}\label{tau6}
\tau_6 (n)=\left( n^2 \right)^2 +2nO(n) +n^2 O(n) =O(n^4 ) .
\end{equation}

Here we will present a more efficient algorithm for obtaining random Sudoku matrix, based on the propositions and algorithms which are examined in the previous sections of this paper. The main point is to obtain $n^2$ in number random $\Pi_n$ matrices (Algorithm \ref{alg4}). For every $\Pi_n$ matrix which is obtained, it has to be checked if it is disjoint with each of the above obtained matrices.The criteria, described in Theorem \ref{th1} and Proposition \ref{disj} are used in the verification. If the obtained matrix is not disjoint with at least one of the above mentioned matrices, it has to be replaced with another randomly generated $\Pi_n$ matrix.

\begin{algorithm}\label{Sudokuend}
Obtaining random Sudoku matrix\\
1) We declare $n^2 \times n^2$ matrix $S$ and initialize all of its elements with 0;\\
2) \textbf{for} $k=1,2,\ldots ,n^2$ \textbf{do}

\textbf{begin}\\
3) We declare $n^2 \times n^2$ matrix $B$ and initialize all of its elements with 0;\\
4) $z:=\mathbf{true}$;\\
5) \textbf{while} $z$ \textbf{do}

\hspace{1cm} \textbf{begin}\\
6) \hspace{1cm} We obtain random matrix $P_k \in \Pi_n$  /* Algorithm \ref{alg4} */ \\
7)\hspace{1cm} We obtain $\Sigma_{n^2}$  matrix $A=\varphi (P_k )$; /* Algorithm \ref{alg5}, where $\varphi$ is defined in Theorem \ref{th1} bijective map */ \\
8) \hspace{1cm} $C:=B+A$;\\
9) \hspace{1cm} \textbf{if} All elements of matrix $C$ are equal to 0 or 1 \textbf{then}
 $P_k =\varphi^{-1} (A)$ is disjoint with each of the matrices $P_1 ,P_2 ,\ldots ,P_{k-1}$;  /* (Theorem \ref{th1})  */

 \hspace{2cm} \textbf{begin}\\
10) \hspace{2cm} $B:=C$;\\
11) \hspace{2cm} $S:=S+k\cdot A$; /* Proposition \ref{disj} */ \\
12) \hspace{2cm} $z:=\mathbf{false}$;

\hspace{2cm} \textbf{end};

\hspace{1cm} \textbf{end};

\textbf{end}.
\end{algorithm}

\end{document}